\documentclass[12pt,bezier]{article}
\usepackage{amsmath,amssymb,amsfonts,euscript,graphicx}

\newtheorem{prethm}{{\bf Theorem}}

\newenvironment{thm}{\begin{prethm}{\hspace{-0.5
               em}{\bf.}}}{\end{prethm}}

\newtheorem{prepro}[prethm]{{\bf Theorem}}

\newtheorem{preprop}[prethm]{{\bf Proposition}}

\newtheorem{precor}[prethm]{{\bf Corollary}}

\newtheorem{preconj}[prethm]{{\bf Conjecture}}

\newtheorem{preremark}[prethm]{{\bf Remark}}

\newtheorem{preexample}[prethm]{{\bf Example}}

\newtheorem{prelem}[prethm]{{\bf Lemma}}

\newtheorem{prelam}{{\bf Lemma}}

\newtheorem{preproof}{{\bf Proof.}}

\title{\large \bf Erratum to: Matching Subspaces in a Field Extension}

\date{}
\begin{document}
%{\large \bf To matching subspaces in a field extension}

\author{{\normalsize{\sc S. Akbari${}^{\mathsf{a}, \mathsf{b}}$}, {\sc M. Aliabadi${}^{\mathsf{a}}$}}
\\
\vspace{3mm}
\\{\footnotesize{${}^{\mathsf{a}}$\it Department of Mathematical Sciences, Sharif University of Technology, }}
{\footnotesize{}}\\{\footnotesize{${}^{\mathsf{b}}$\it School of Mathematics, Institute for Research in Fundamental Sciences,}}\\{\footnotesize{\it P.O. Box \rm{19395-5746}, Tehran, Iran}}\\
%\\
{\footnotesize{}}\\
{\footnotesize{$\mathsf{s\_akbari@sharif.edu}$\hspace*{1cm}
$\mathsf{mohsenmath88@gmail.com}$}}\\
}

\maketitle

\begin{abstract}
{
Let $K\subset L$ be a field extension. The field $L$ is said to have the linear matching property if, for every $n\geq1$ and every $n$-dimensional $K$-subspaces $A$ and $B$ of $L$ with $1\not\in B$, the subspace $A$ is matched to $B$.
Recently, it has been proved that if $K\subset L$ be a field extension, then $L$ has linear matching property if and only if $L$ is purely transcendental or is an extension of $K$ of prime degree. In this note we provide a counterexample for this result.

}
\end{abstract}


\maketitle
\section{Introduction}
Throughout the paper, we shall say \textit{field} for a division ring. 
If $K\subset L$ is a field extension, then we say that $K$ is \textit{central} in $L$ if for all $\lambda\in K$, $x\in L$\quad $\lambda x=x\lambda$.
Suppose that $A$ and $B$ are two $n$-dimensional $K$-subspaces of $L$. Let $\mathcal A=\{a_1,\ldots,a_n\}$ and $\mathcal B=\{b_1,\ldots,b_n\}$ be two bases of $A$ and $B$, respectively. We say that $\mathcal A$ is \textit{matched} to $\mathcal B$ if
\[
a_ib\in A~~\Rightarrow~~b\in\langle b_1,\ldots,\hat{b_i},\ldots,b_n\rangle
\]
for all $b\in B$ and $i=1,\ldots,n$, where $\langle b_1,\ldots,\hat{b_i},\ldots,b_n\rangle$ is the subspace of $B$ spanned by the set $\mathcal B\setminus \{b_i\}$; equivalently, if 
\[
a_i^{-1}A\cap B\subset\langle b_1,\ldots,\hat{b_i},\ldots,b_n\rangle
\]
for all $i=1,\ldots,n$.\\
Let $A$ and $B$ be two $n$-dimensional $K$-subspaces in $L$. We say that $A$ is \textit{matched} to $B$ if every basis $\mathcal{A}$ of $A$ can be matched to a basis $\mathcal{B}$ of $B$.
Let $K\subset L$ be a field extension. Then $L$ is said to have the \textit{linear matching property} if, for every $n\geq1$ and every $n$-dimensional $K$-subspaces $A$ and $B$ of $L$ with $1\not\in B$, the subspace $A$ is matched to $B$.
If $K\varsubsetneq L$ be a field extension, where $K$ is commutative and central in $L$, then $n_0(K,L)$ denotes the smallest degree of an intermediate field extension $M$, $K\varsubsetneq M\subset L$.
\\
The following theorem was proved in [1]:

\begin{thm}
Let $K\subset L$ be a field extension, with $K$ commutative and central in $L$. Then $L$ has linear matching property if and only if $L$ is purely transcendental or is an extension of $K$ of prime degree.
 \end{thm}
There is a gap in the proof of Theorem 1. In Line $2$ of the proof, the authors claimed that if $K\subset L$ is a field extension and $[ L : K]$ is not prime, then there exists $a\in L$ such that $K\varsubsetneq K(a)\varsubsetneq L$. But this claim does not hold. To see this we note that there is a field extension $Q\subset L$ of degree 4 with no proper intermediate subfield, see [2].
By the definition $n_0(Q,L)=4$. Let $A$ and $B$ be two $n$-dimensional $K$-subspaces in $L$ such that $1\not\in B$. Since $1\not\in B$, we have $n < 4$. Thus $n < n_0(Q,L)$. Now, Theorem $5.3$ of [1] implies that $A$ is matched with $B$. This contradicts the statement of Theorem 1.\\
{}


\begin{thebibliography}{}{\small

\bibitem{} S. Eliahou, C. Lecouvey, Matching subspaces in a field extension, J. Algebra 324 (2010) 3420-3430.

\bibitem{} H. Marksaitis, Some remarks on subfields of algebraic number fields, Lithuanian Mathematical Journal, Vol. 35. No. 2. 1995.

}\end{thebibliography}
\end{document}